\documentclass{endm}
\usepackage{endmmacro}
\usepackage{graphicx,enumerate}

\newcommand\cB{\mathcal B}

\newcommand\GF{\operatorname{GF}}

\begin{document}

\begin{verbatim}\end{verbatim}\vspace{2.5cm}

\begin{frontmatter}

\title{Negative Circles in Signed Graphs: \\[6pt]A Problem Collection}

\author{Thomas Zaslavsky\thanksref{myemail}}
\address{Department of Mathematical Sciences\\ Binghamton University (State University of New York)\\ Binghamton, NY 13902-6000, U.S.A.}
\thanks[myemail]{Email:
   \href{mailto:zaslav@math.binghamton.edu} {\texttt{\normalshape
   zaslav@math.binghamton.edu}}} 

\begin{abstract}
I propose that most problems about circles (cycles, circuits) in ordinary graphs that have odd or even length find their proper setting in the theory of signed graphs, where each edge has a sign, $+$ or $-$.  Even-circle and odd-circle problems correspond to questions about positive and negative circles in signed graphs.  (The sign of a circle is the product of its edge signs.)  I outline questions about circles in signed graphs, that seem natural and potentially important.
\end{abstract}


\end{frontmatter}

\section{Beginning}\label{beg}

Problems about circles (cycles, circuits) in ordinary graphs have attracted much attention over the years.  
Problems about circles of odd length, or even length, form a small but increasing important part of that area.  I propose that the proper setting for questions involving parity of circles is, mostly, the theory of signed graphs, where each edge has a sign.  Even-circle and odd-circle problems for unsigned graphs generalize, respectively, to questions about positive and negative circles in signed graphs.  
Here I outline a framework for questions about circles in signed graphs that seem natural and potentially valuable.  Some of the questions arise from other research; most are simply basic structural questions needed for a better understanding of signed graphs.

A \emph{signed graph} is $\Sigma = (\Gamma,\sigma)$, where $\Gamma$ is a graph and $\sigma: E(\Gamma) \to \{+,-\}$, the \emph{signature}, is a sign function on the edges.  
I assume all graphs are simple and I write $n := |V|$, the order of the graph.
The \emph{sign of a circle} is the product of the signs of its edges. 
A crucial property of a signed graph is the list of its positive circles.  Some elementary properties that depend only on that list:

\emph{Balance}:  $\Sigma$ is balanced if all circles are positive.  Otherwise it is unbalanced, or (in physics) frustrated.  For instance, the all-positive signed graph $(\Gamma,+)$ is balanced.

\emph{Balancing edge}: $\Sigma$ is unbalanced but $\Sigma\setminus e$ is balanced.

\emph{Frustration}: The \emph{frustration index} $l(\Sigma)$ is the minimum number of edges whose deletion makes $\Sigma$ balanced; that is, eliminates all negative circles.
Similarly, the \emph{frustration number} $l_0(\Sigma)$ is the minimum number of vertices whose deletion makes $\Sigma$ balanced.

Suppose $\Sigma$ is all negative, $\Sigma=(\Gamma,-)$.  
Then the negative circles are the odd circles and $\Sigma$ is balanced if and only if $\Gamma$ is bipartite.  

Answers to the first six questions are known; but most of the problems are unsolved.  If there is no ``\emph{Ans.}''\ line, the answer is unknown (at least to me).  Each conjecture is an educated guess and may well be mistaken.

\section{The Structure of the Class of Negative Circles}\label{N}

The solution to Problem 1 is essential to all work of this kind.  The solution to Problem 2 is probably not needed at all, but it shows the difference in complexity between circles and chordless circles.  Any question about circles can have the word ``chordless'' added to make a new question that is worth working on---but harder.

\begin{enumerate}[1.]

\item  Can a given set $\cB$ of circles in $\Gamma$ be the negative circles of a signature?
\\\emph{Ans.}  Easy:  $\cB$ is a negative circle set iff every theta subgraph of $\Gamma$ has, among its three circles, either one or three that belong to $\cB$ \cite{CSG}.

\item  Can a given set of chordless circles in $\Gamma$ be the negative chordless circles of a signature?
\\\emph{Ans.}  There are infinitely many forbidden subgraphs of a finite number of kinds.  This important theorem is due to Truemper \cite{Tru}.

\end{enumerate}

\section{Edges and Vertices in Circles}

\subsection{Edges in Negative Circles}\label{ENC}

\begin{enumerate}[E1.]

\item  In $\Sigma$, does a given edge $e$ belong to a negative circle?
\\\emph{Ans.}  Easy:  Yes, if and only if $e$ is in an unbalanced block (Harary \cite{LB}).
\label{Eneg}

\item  Assume $\Sigma$ is 2-connected and unbalanced.  Is every pair of edges in a common negative circle?
\\\emph{Ans.}  Not every pair is; it depends on $\Sigma$.  How does it depend?  For which $\Sigma$ is it true?  This is one of the most basic and important questions.  I believe the answer is that every pair is, if $\Sigma$ is 3-connected, and otherwise it depends on the Tutte 3-decomposition of $\Gamma$ \cite{Tbook}.  Possibly, every pair is, if $\Sigma$ does not have a balancing edge.
\label{Epairneg}

\item  In $\Sigma$, is a certain edge $e$ in a unique negative circle?
\\\emph{Ans.}  Medium hard:   It depends on the Tutte 3-decomposition of $\Gamma$ into 2-connected subgraphs, and how those subgraphs are signed (Behr \cite{Behr}).
\label{E!neg}

\item  In $\Sigma$, find the set of all edges $e$ such that $e$ belongs to a unique negative circle.
\\\emph{Ans.}  Similar to E\ref{E!neg} (Behr \cite{Behr}).
\label{Eall!neg}

\item  In a connected $\Sigma$, which edges $e$ belong only to negative circles?
\\\emph{Ans.}  Unknown.  \emph{Conjecture}:  They are the isthmi and the balancing edges.
\label{Enopos}

\label{Elast}
\end{enumerate}

For the chordless-circle version of E\ref{Eneg} there is a simple algorithm by Marinelli and Parente \cite{MarPar} (and they ask for a better algorithm) but no answer in terms of graph structure.

\subsection{Edges in Positive Circles}  \label{EPC}

To me, positive circles seem less fundamental than negative circles---possibly because positive sign seems like the default sign.  Still, they can be as important and since they obey a different rule, their properties are not the same.
I particularly want the answers to EP\ref{EPpos}--\ref{EPall!pos} because they will illuminate the difference between positive and negative circles.

\begin{enumerate}[EP1.]

\item  In $\Sigma$, is a certain edge $e$ in a positive circle?
\\\emph{Ans.}  Not known.  This question is the negation of E\ref{Enopos}.
\label{EPpos}

\item  Assume $\Gamma$ is 2-connected.  In $\Sigma$, does every pair of edges lie in a positive circle?
\\\emph{Ans.}  This should be similar to EP\ref{Epairneg}.
\label{EPpairpos}

\item  In $\Sigma$, is a certain edge $e$ in a unique positive circle?
\\\emph{Ans.}  It should broadly resemble E\ref{E!neg} but have distinctive features.
\label{EP!pos}

\item  In $\Sigma$, find the set of all edges $e$ such that $e$ belongs to a unique positive circle.
\\\emph{Ans.}  Similar to EP\ref{EP!pos}.
\label{EPall!pos}

\item  In $\Sigma$, which edges $e$ belong only to positive circles?
\\\emph{Ans.}  Easy:  The edges in balanced blocks.  This is the negation of E\ref{Eneg}.
\label{EPnoneg}

\end{enumerate}

\subsection{Vertices in Negative Circles}\label{VNC}

\begin{enumerate}[V1.]

\item  In $\Sigma$, is a certain vertex $v$ in a negative circle?
\\\emph{Ans.}  Yes, if and only if $v$ is in an unbalanced block.
\label{Vneg}

\item  In $\Sigma$, does a certain pair of vertices, $v, w$, lie in a negative circle?
\\\emph{Ans.}  Yes, if and only if they are in an unbalanced block together.
\label{Vpairneg}

\item  In $\Sigma$, is a certain vertex $v$ in a unique negative circle?
\\\emph{Ans.}  Similar to E\ref{E!neg}.  Solved by Behr \cite{Behr}.
\label{V!neg}

\item  In $\Sigma$, does vertex $v$ belong only to negative circles?
\\\emph{Ans.}  
\emph{Conjecture}: In every unbalanced block $B$ that contains it, $v$ must be divalent and a balancing vertex (that is, $B \setminus v$ is balanced).%
\label{Vnopos}

\label{Vlast}
\end{enumerate}

On the whole, though not always, vertex answers should follow from edge answers.

\subsection{Vertices in Positive Circles}\label{VPC}

Questions VP\ref{Vneg}--\ref{Vlast} are V\ref{Vneg}--\ref{Vlast} with the obvious change of ``negative'' to ``positive''.  The answers may also be similar, but VP\ref{Vneg}--\ref{Vpairneg} are not as easy.

VP\ref{Vnopos} is easy: $v$ must belong only to balanced blocks.  That follows from the fact that, in an unbalanced block, every edge belongs to a negative circle (Harary; see E\ref{Eneg}).

\section{Packing and Covering}

\subsection{Packing Circles}\label{P}

Work on odd or even circles in unsigned graphs ought to generalize.

\begin{enumerate}[P1--9.]

\item[P1.]  What is the maximum number of pairwise disjoint negative circles in $\Sigma$?
\\\emph{Ans.}  Unknown.  It is obvious that this number is $\leq l_0(\Sigma)$.  Which signed graphs have equality?  
Slilaty \cite{Sli} treats two circles; Hochst\"attler et al.\ \cite{HNP} have an algorithm.  
Some parity papers:  Berge and Reid \cite{BR} treat odd circles; Kr\'al' et al.\ \cite{Kral} treat odd circles in planar graphs.  
\label{Pnegnumber}

\item[P2.]  Find a maximum set of pairwise disjoint negative circles.
\label{Pnegset}

\item[P3--4.]  The same as P1--2, for positive circles.
\label{Pposset}

\item[P5--8.]  The same as P1--4 but for edge-disjoint circles.
\\\emph{Ans.}  Unknown.  The maximum number of negative circles obviously is $\leq l(\Sigma)$.  When is there equality?  

\label{Plast}
\end{enumerate}

\subsection{Covering by Circles}\label{C}

\begin{enumerate}[C1--4.]

\item[C1--4.]  Like P\ref{Pnegnumber}--\ref{Pposset} but for the minimum number, or minimum sets, of negative (or positive) circles that cover all the vertices of $\Sigma$.

\item[C5--8.]  Like C1--4, for circles that cover the edges of $\Sigma$.

\item[C9.]  Are there duality relations between packing and covering numbers?

\label{Clast}
\end{enumerate}

\subsection{Decomposition into Circles}\label{D}

These are suggested by the theorem that a connected graph decomposes into circles iff it is Eulerian.  (\emph{Decomposing} a graph means partitioning its edge set.)  D\ref{Dnegdecomp}--\ref{Dposdecomp} seem very hard.  D\ref{Dnegptn} is open-ended.

\begin{enumerate}[D1.]

\item  Can $\Sigma$ be decomposed into negative circles?
\label{Dnegdecomp}

\item  Can $\Sigma$ be decomposed into positive circles?
\label{Dposdecomp}

\item  Is there an interesting Euler-type property of a connected signed graph related to either D\ref{Dnegdecomp} or D\ref{Dposdecomp}?
\label{Dnegptn}

\label{Dlast}
\end{enumerate}

\section{Counting Negative Circles}\label{CN}

The \emph{negative circle vector} is $c^-(\Sigma) = (c_3^-,c_4^-,\ldots,c_n^-) \in {\mathbb R}^{n-2}$, where  $c_l^-(\Sigma)$ is the number of negative circles of length $l$.

\begin{enumerate}[CN1.]

\item  Characterize the set of numbers of negative circles of some fixed length of all signatures of $K_n$; that is, $\{c_l^-(K_n,\sigma) : \sigma$ is a signature of  $K_n\}$ for some fixed $l$, $3\leq l \leq n$.
\\\emph{Ans.}  Very recently there are remarkably strong results on the possible numbers of negative triangles  (Kittipassorn and M\'esz\'aros \cite{KittiMesz}).  I am not aware of any results about longer circles.

\item  Characterize the possible numbers of negative circles of a fixed length of all signatures of $K_{r,s}$.

\item  Characterize the set ${\mathbf C}$ of negative circle vectors of all signatures of a fixed graph $\Gamma$.
\\\emph{Ans.}  There are partial results for complete graphs $K_n$ by Popescu and Tomescu, e.g.\ in  \cite{Pop1,Pop2,Tomescu}.  
Schaefer \cite{S} finds that $\dim {\mathbf C}$ is the largest it could possibly be for $\Gamma = K_n$ and $K_{r,s}$:  that is, $n-2$ and $\min(r,s)-1$, respectively.  

\label{Nlast}
\end{enumerate}

\section{Structure of the Signed Graph}\label{SSG}

\begin{enumerate}[S1.]

\item  Assume $\Sigma$ has a Hamiltonian circle and is unbalanced.  Is there a negative Hamiltonian circle?  A positive one?
\\\emph{Ans.}  Unknown.  \emph{Conjecture}:  Most $\Sigma$ with a Hamiltonian circle have both.  The exceptions are unbalanced necklaces (circles) of balanced blocks \cite[Theorem 5.9]{SG}, and have only negative Hamiltonian circles.

\item  Is there a positive, or negative, circle $C$ such that $\Sigma \setminus E(C)$ is disconnected, or separable, or 2-separable, or 2-connected?
\\\emph{Ans.}  Conlon \cite{Conlon} proved that if $\Gamma$ is 3-connected, there is an even circle $C$ such that $\Gamma \setminus E(C)$ is 2-connected.  Fujita and Kawarabayashi \cite{FuKawara} have a similar theorem for $\Gamma \setminus V(C)$.
Do these generalize to signed graphs?  What definition of connectivity of a signed graph is suggested?

\item  What are the bridges (in the sense of Tutte) of a negative or positive circle?  For instance, does the circle have many chords?
\\\emph{Ans.}  Voss \cite{Voss} studied chords and other properties of circles in $\Gamma$ of given parity.  Which of these generalize to circles of given sign in $\Sigma$?

\label{Slast}
\end{enumerate}


\end{document}